\documentclass[english]{article}
\usepackage[T1]{fontenc}
\usepackage[latin9]{inputenc}
\usepackage{units}
\usepackage{amsmath}
\usepackage{amssymb}

\makeatletter
\newcommand{\lyxaddress}[1]{
\par {\raggedright #1
\vspace{1.4em}
\noindent\par}
}
\newenvironment{lyxcode}
{\par\begin{list}{}{
\setlength{\rightmargin}{\leftmargin}
\setlength{\listparindent}{0pt}
\raggedright
\setlength{\itemsep}{0pt}
\setlength{\parsep}{0pt}
\normalfont\ttfamily}%
 \item[]}
{\end{list}}

\makeatother

\usepackage{babel}
\begin{document}

\title{Dimensional Analysis: A Centenary Update}

\author{Dan Jonsson}
\maketitle
\begin{abstract}
It is time to renew old ways of thinking about dimensional analysis.
Specifically,\emph{ }more than $n-r$ invariants and more than one
functional relation between invariants need to be considered simultaneously.
Thus generalized, dimensional analysis can yield more information
than previously recognized.
\end{abstract}
Buckingham's $\Pi$ theorem \cite{bucking} and its predecessors \cite{vaschy,riab}
were formulated 100 years or more ago. The basic principles of dimensional
analysis have remained unchanged since then. Yet, a careful investigation
reveals that dimensional analysis rests on presuppositions which unnecessarily
limit the scope and power of the analysis. This is implied by the
examination of basic notions related to dimensional analysis in \cite{dan};
this article will explain more explicitly how dimensional analysis
should be extended to transcend its self-imposed limitations.

\section{Dimensional analysis should be liberated from some traditional constraints}

Dimensional analysis as depicted in well-known classical expositions
\cite{bridgma,sedov,baren}, can be summarized as follows. We want
to express a dependent quantity $q$ as a function of independent
quantities $q_{1},\ldots,q_{n-1}$, where $n>1$; thus, we assume
that $q=\Psi\!\left(q_{1},\ldots,q_{n-1}\right)$, or equivalently
and sometimes more conveniently, 
\begin{equation}
q_{1}=\Psi\!\left(q_{2},\ldots,q_{n}\right).\label{eq:Psi-rel}
\end{equation}
We also assume that the quantities $q_{1},\ldots,q_{n}$ can all be
expressed in terms of one or more fundamental units of measurement,
corresponding to $m$ dimensions such as Length, Time and Mass, where
$1\leqq m<n$. Dimensional analysis makes it possible to find $n-r$
so-called $\pi$-\emph{groups} $\pi_{1},\ldots,\pi_{n-r}$ such that
(\ref{eq:Psi-rel}) can be written as 
\begin{equation}
\pi_{1}=\Phi\!\left(\pi_{2},\ldots,\pi_{n-r}\right),\label{eq:Phi-rel}
\end{equation}
for some $r$ such that $0\leq r\leq m$. The $\pi-$groups are invariant
products of powers of $q_{1},\ldots,q_{n}$, meaning that the numerical
values of these products do not depend on the fundamental units of
measurement used to express $q_{1},\ldots,q_{n}$. Relations of the
form (\ref{eq:Phi-rel}) can be rewritten as
\[
q_{1}^{c}=\prod_{j=1}^{r}q_{i_{j}}^{c_{j}}\Phi\!\left(\pi_{1},\ldots,\pi_{n-r}\right)\quad\mathrm{or}\quad q_{1}=\prod_{j=1}^{r}q_{i_{j}}^{c_{j}/c}\Phi\!\left(\pi_{1},\ldots,\pi_{n-r}\right)^{1/c},
\]
where $c,c_{j}$ are integers and $q_{1},q_{i_{j}}$ are distinct
quantities.

In many contemporary expositions of dimensional analysis, the way
in which $q_{1},\ldots,q_{n}$ are expressed in terms of fundamental
units of measurement or corresponding dimensions is summarized in
a \emph{dimensional matrix} $\left[a_{ij}\right]$, where $a_{ij}$
is the dimensionality of $q_{j}$ relative to the dimension $D_{i}$,
meaning that $\left[q_{j}\right]=D_{1}^{a_{1j}}\cdots D_{i}^{a_{ij}}\cdots D_{m}^{c_{mj}}$
\cite{dan}. The invariants ($\pi-$groups) are obtained from this
matrix. It is not difficult to show that $r$ is equal to the rank
of the dimensional matrix. The introduction of notions and techniques
from linear algebra does not change the general way of thinking about
dimensional analysis, however.

Dimensional analysis as described above has two limitations:

(P1) \quad{}\emph{Not more than one relation of the form (\ref{eq:Phi-rel})
is considered}.

(P2) \quad{}\emph{Not more than $n-r$ invariants are considered}.

The main purpose of this article is to present a generalized form
of dimensional analysis which is not constrained by these limitations.
Before proceeding, it should be noted, though, that there exists a
tension between the traditional, abstract description of dimensional
analysis and the method used in practice. This method revolves around
the notion of \emph{repeating variables} (specifically, repeating
quantities), and is based on the observation that every invariant
in a relation of the form (\ref{eq:Phi-rel}) can be written as 
\[
\frac{p^{c}}{p_{1}^{c_{1}}\cdots p_{r}^{c_{r}}}\qquad\mathrm{or}\qquad\frac{p}{p_{1}^{c_{1}/c}\cdots p_{r}^{c_{r}/c}},
\]
where $p,p_{i}$ are distinct quantities, $c,c_{i}$ are integers.
Here, each $p$ which occurs in a numerator occurs in exactly one
invariant in the relation, while each $p_{i}$ which occurs in a denominator
may occur in all invariants in the relation -- hence, repeating variables
-- but \emph{the set of quantities} \emph{occurring in the denominators
can in general be chosen in more than one way}. 

As we shall see, different sets of repeating variables give different
relations of the form (\ref{eq:Phi-rel}). Hence, there may be more
than one such relation, and as the $n-r$ invariants in one relation
are not the same as the $n-r$ invariants in another relation, the
total number of invariants may be greater than $n-r$. The emphasized
fact is thus an anomaly from the point of view of dimensional analysis
in the tradition of the $\Pi$ theorem. Simplifying and idealizing,
one can say that, apart from quietly ignoring the anomaly, three ways
of dealing with it are found in the literature:
\begin{enumerate}
\item Some authors try to re-establish uniqueness by suggesting criteria
for choosing the 'right' repeating variables, although they may not
insist that finding one 'right' set of repeating variables is always
possible.
\item Some authors acknowledge non-uniqueness but downplay it, arguing that
essentially the same result is obtained no matter which formally well-behaved
set of repeating variables is chosen.
\item Some authors accept non-uniqueness in practice, thus effectively abandoning
(P1) and (P2), but pay lip service to the traditional formulation
of dimensional analysis, thereby disconnecting theory and practice.
\end{enumerate}
These positions will be described more fully below, relating each
one to the approach proposed here.

\section{\label{sub:1-1}Invariants versus active invariants}

An example, adapted from Sedov \cite{sedov}, will be used to introduce
the new way of thinking about dimensional analysis. Consider a fluid
flowing through a cylindrical pipe; let $\nicefrac{\triangle P}{\ell}$
be the pressure drop per unit length, $\rho$ the density of the fluid,
$\mu$ the viscosity of the fluid, $d$ the diameter of the pipe,
and $u$ the (mean) velocity of the fluid. The dimensional matrix
for these quantities is
\[
\begin{array}{cc}
\begin{array}{c}
\\
L\\
T\\
M
\end{array} & \left\Vert \begin{array}{ccccc}
\nicefrac{\triangle P}{\ell} & \rho & \mu & d & u\\
-2 & -3 & -1 & 1 & 1\\
-2 & 0 & -1 & 0 & -1\\
1 & 1 & 1 & 0 & 0
\end{array}\right\Vert .\end{array}
\]
The rank of this matrix is $3$. There is a well-known invariant involving
the quantities $\rho,\mu,d,u$, namely the Reynolds number
\[
\frac{\rho du}{\mu}=\mathsf{Re}=\Pi_{1}.
\]
 Sedov also considers another invariant
\[
\frac{\left(\nicefrac{\triangle P}{\ell}\right)d}{\rho u^{2}}=\Pi.
\]
and we now have $n-r=5-3=2$ invariants. Sedov asserts that $\mathsf{Re}$
is the only possible invariant involving $\rho,\mu,d$ and $u$, which
means that any invariant which depends on $\rho,\mu,d,u$ actually
depends on $\rho du\mu^{-1}$. Thus, assuming that $\nicefrac{\triangle P}{\ell}=\Psi\!\left(\rho,\mu,d,u\right)$,
we have $\Pi=\Phi\!\left(\Pi_{1}\right)$, or
\[
\nicefrac{\triangle P}{\ell}=\rho u^{2}d^{-1}\Phi\!\left(\mathsf{Re}\right).
\]

However, we have not looked for invariants systematically. We stopped
after finding two invariants, guided by the traditional way of thinking
about dimensional analysis. But there is actually no ground for assuming
that these are the only invariants for this dimensional matrix, or
even the only invariants of interest, so let us see what happens if
we look for more invariants.

In a sense, one should not look for all invariants, however. To begin
with, if $\pi$ is an invariant then $\pi^{k}$ is also an invariant
for any non-zero integer $k$; yet, all these invariants are essentially
the same. To better understand what this means, consider a product
of powers of quantities $q_{i_{1}}^{c_{1}}\ldots q_{i_{k}}^{c_{k}}$,
and let $\mathbf{v}_{j}$ be the (column) vector corresponding to
$q_{i_{j}}$ in the dimensional matrix. As explained in Section 3,
$q_{i_{1}}^{c_{1}}\ldots q_{i_{k}}^{c_{k}}$ is an invariant if and
only if the equation with integer coefficients
\[
c_{1}\mathbf{v_{1}}+\ldots+c_{k}\mathbf{v}_{k}=\mathbf{0}
\]
holds. It is clear that this equation still holds, and still has integer
coefficients, if the integers $c_{1},...,c_{k}$ are multiplied by
a non-zero integer, or divided by a common divisor of $c_{1},\ldots,c_{k}$.
Thus, we can represent all such equations by an equation where $c_{1},\ldots,c_{k}$
are relatively prime. There are exactly two such equations, which
can be obtained from each other through multiplying $c_{1},...,c_{k}$
by $-1$. The corresponding set of invariants can then be represented
by \emph{pairs of minimal invariants}, each of which can be obtained
from the other by a sign flip. Expressed otherwise, there are sets
of equivalent invariants such that each set can be represented --
uniquely up to a sign flip -- by an invariant $q_{i_{1}}^{c_{1}}\ldots q_{i_{k}}^{c_{k}}$
such that $c_{1},\ldots,c_{k}$ are relatively prime.

The two invariants considered by Sedov could thus be regarded as a
set of pairs of minimal invariants
\[
\left\{ \left(\frac{\rho du}{\mu}\right)^{\pm1},\left(\frac{\left(\nicefrac{\triangle P}{\ell}\right)d}{\rho u^{2}}\right)^{\pm1}\right\} .
\]
 In simplified notation, this set can be denoted
\[
\left\{ \frac{\rho du}{\mu},\frac{\left(\nicefrac{\triangle P}{\ell}\right)d}{\rho u^{2}}\right\} ^{\pm1}.
\]

So far, we have only succeeded in reformulating the original question
\emph{''Are there more than $n-r$ invariants?''} into the more
sophisticated question \emph{``Are there more than $n-r$ pairs of
minimal invariants?''}, however. The situation will be clarified
in Section 3; a more informal discussion will suffice here. 

Recall that the dimensional matrix given has rank $3$, and let us
accept that the sets of possible repeating quantities that we need
to consider are $\left\{ \rho,\mu,d\right\} $, $\left\{ \rho,\mu,u\right\} $,
$\left\{ \rho,d,u\right\} $ and $\left\{ \mu,d,u\right\} $. The
corresponding homogeneous systems of linear equations have minimal
integer solutions defining the following invariants:

\[
\left\{ \frac{\nicefrac{\triangle P}{\ell}}{\rho^{-1}\mu^{2}d^{-3}},\frac{\nicefrac{\triangle P}{\ell}}{\rho^{2}\mu^{-1}u^{3}},\frac{\nicefrac{\triangle P}{\ell}}{\rho d^{-1}u^{2}},\frac{\nicefrac{\triangle P}{\ell}}{\mu d^{-2}u},\frac{\rho ud}{\mu}\right\} ^{\pm1}.
\]
We can choose one invariant in each pair without loss of generality.
Accordingly, dimensional analysis produces four -- not one -- possible
ways of writing the relation $\nicefrac{\triangle P}{\ell}=\Psi\!\left(\rho,\mu,d,u\right)$,
and we can take these representations to be:

\[
\begin{cases}
\nicefrac{\triangle P}{\ell}=\rho^{-1}\mu^{2}d^{-3}\,\phi_{1}\!\left(\rho ud\mu^{-1}\right) & (\mathfrak{a})\\
\nicefrac{\triangle P}{\ell}=\rho^{2}\mu^{-1}u^{3}\,\phi_{2}\!\left(\rho ud\mu^{-1}\right) & (\mathfrak{b})\\
\nicefrac{\triangle P}{\ell}=\rho d^{-1}u^{2}\,\phi_{3}\!\left(\rho ud\mu^{-1}\right) & (\mathfrak{c})\\
\nicefrac{\triangle P}{\ell}=\mu d^{-2}u\,\phi_{4}\!\left(\rho ud\mu^{-1}\right) & (\mathfrak{d})
\end{cases}.
\]
As we have seen, Sedov gives formula $\left(\mathfrak{c}\right)$,
and other authors seem to have followed in his footsteps \cite{baren}. 

When the internal roughness of the pipe can be disregarded as assumed
above, $\left(\mathfrak{c}\right)$ is equivalent to the Darcy-Weisbach
equation, which can be written in notation similar to that used above
as 
\[
\triangle P=f\frac{\ell}{d}\frac{\rho u^{2}}{2}.
\]
Here, $f=f\!\left(\frac{\rho ud}{\mu},\frac{\epsilon}{d}\right)$
is a ``dimensionless'' quantity and $\epsilon$ is a quantity of
dimension $L$ expressing the roughness of the pipe, with 
\[
f\!\left(\frac{\rho ud}{\mu},\frac{\epsilon}{d}\right)=2\Phi_{3}\!\left(\frac{\rho ud}{\mu}\right)
\]
for the idealized case $\frac{\epsilon}{d}=0$.

Is ($\mathfrak{c}$) the 'right' formula, then? Can we reject the
three other formulas? No, no formula is 'wrong'. All four formulas
are equally correct and in fact interchangeable for the somewhat surprising
reason that all contain less information than they seem to do at first
sight. A closer look at formula ($\mathfrak{c}$), for example, reveals
that it is not possible to conclude that $\nicefrac{\triangle P}{\ell}$
is proportional to $\rho$ and $u^{2}$ and inversely proportional
to $d$, because $\rho,$ $\mu$ and $d$ also appear as arguments
of $\Phi_{3}$. Corresponding conclusions apply to formulas ($\mathfrak{a}$),
($\mathfrak{b}$) and ($\mathfrak{d}$). The four formulas are different
because the effect of \emph{different} quantities are ``hidden inside''
the functions $\Phi_{1}$ through $\Phi_{4}$. Specifically, the effects
of $u$, $d$, $\mu$ and $\rho$ are hidden in formulas ($\mathfrak{a}$),
($\mathfrak{b}$), ($\mathfrak{c}$) and ($\mathfrak{d}$), respectively.
Let us consider two more cases to clarify the situation.

For Reynolds numbers less than a critical value $\mathsf{Re}_{c}$,
the flow through the pipe is laminar, which means that the flow is
non-accelerated. This implies that the flow is not affected by $\rho$
\cite{sedov}. (Simply stated, acceleration depends on inertia, which
depends on mass, which depends on density.) Thus, we obtain the following
dimensional matrix for laminar flow through a pipe:
\[
\begin{array}{c}
\\
L\\
T\\
M
\end{array}\left\Vert \begin{array}{cccc}
\nicefrac{\triangle P}{\ell} & \mu & d & u\\
-2 & -1 & 1 & 1\\
-2 & -1 & 0 & -1\\
1 & 1 & 0 & 0
\end{array}\right\Vert .
\]
In this case, we need to consider only one pair of minimal invariants,
\[
\left(\frac{\nicefrac{\triangle P}{\ell}\, d^{2}}{\mu u}\right)^{\pm1},
\]
and the relation $\nicefrac{\triangle P}{\ell}=\Psi\left(\mu,d,u\right)$
can be written in the form
\begin{equation}
\nicefrac{\triangle P}{\ell}=\frac{\mu u}{d^{2}}\Phi\!\left(\right)=K\frac{\mu u}{d^{2}},\label{eq:lamin}
\end{equation}
 where $K$ is a scalar constant. 

This is obviously a special case of ($\mathfrak{d}$). It is worth
noting that the Hagen--Poiseuille equation, which can be written in
notation similar to that used here as
\[
\triangle P=C\frac{\ell\mu Q}{d^{4}},
\]
where $C$ is a scalar constant and $Q$ the volumetric flow rate,
is equivalent to (\ref{eq:lamin}), since $Q=u\left(\pi/4\right)d^{2}$.
Hence, this well-known equation can be regarded as a special case
of $\left(\mathfrak{d}\right)$.

Using the fact that
\[
K\frac{\mu u}{d^{2}}=\nicefrac{\triangle P}{\ell}=\frac{\rho u^{2}}{d}\,\Phi_{3}\!\left(\frac{\rho ud}{\mu}\right)
\]
for laminar flow, we can determine $\Phi_{3}$ for corresponding values
of $\mathsf{Re}$, and we obtain $\Phi_{3}\!\left(x\right)=K/x$ for
$x<\mathsf{Re}_{c}$. The remaining three functions $\Phi_{i}$ can
be derived in the same way, and we have
\[
\begin{cases}
\Phi_{1}\!\left(x\right)=Kx\\
\Phi_{2}\!\left(x\right)=K/x^{2}\\
\Phi_{3}\!\left(x\right)=K/x\\
\Phi_{4}\!\left(x\right)=K
\end{cases}
\]
for laminar flow.

A case of particular historical interest will also be mentioned. In
his article from 1883, \textquotedbl{}An experimental investigation
of the circumstances which determine whether the motion of water shall
be direct or sinuous, and of the law of resistance in parallel channels\textquotedbl{},
Reynolds \cite{reyn-1} not only discussed the critical Reynolds number
(as it came to be called), but also stated a law concerning the pressure
drop connected with the flow of water through a cylindrical pipe.
In notation similar to that used here, his law is
\[
\nicefrac{\triangle P}{\ell}\,\frac{d^{3}}{\mu_{\theta}^{2}}=\mathrm{\Phi\!}\left(\frac{du}{\mu_{\theta}}\right),
\]
where $\mu_{\theta}$ is the viscosity of water at temperature $\theta$
(p. 973). Considering that the only fluid used by Reynolds was water,
and that the temperature of water affects its viscosity much more
than its density, so that $\rho$ can be regarded as a constant, Reynolds
formula is obviously a special case of\textbf{ }($\mathfrak{a}$). 

We have shown, then, that two of the representations of the relation
$\nicefrac{\triangle P}{\ell}=\Psi\!\left(\rho,\mu,d,u\right)$ correspond
to well-known formulas in fluid dynamics, while one representation,
although known to Reynolds, seems to have been forgotten, and one
representation seems to have escaped notice altogether.

To return to the main theme, we conclude that in the interest of clarity
a distinction should be made between \emph{invariants} and \emph{active}
\emph{invariants}. For example, in the main example we consider five
invariants -- or pairs of minimal invariants -- but there are only
two active invariants in each of the four relations derived. It is
clear why this distinction is obscured in traditional dimensional
analysis -- it is not relevant when only one relation is considered.

\subsection*{Remark}

Some textbooks treating dimensional analysis suggest criteria for
choosing repeating quantities in order to help students to choose
the 'right' set of repeating quantities, usually mixing 'formal' and
'substantial' criteria \cite{cengel}. 'Formal' criteria are criteria
such as ``Do not assign two quantities with the same dimensions to
the same set of repeating quantities''. Such criteria ensure that
the mathematical assumptions underlying dimensional analysis are satisfied;
specifically, they guarantee that a set of repeating quantities is
a maximal set of independent quantities and does not contain the dependent
quantity. 'Substantial' criteria are criteria such as ``If possible,
choose a simple quantity such as a length or a time as a repeating
quantity''. These criteria are meant to help students choose the
'right' groups of repeating quantities among the groups which satisfy
the 'formal' criteria, so that, for example, one of the relations
$\left(\mathfrak{a}\right)$ -- $\left(\mathfrak{d}\right)$ is designated
as the 'right' one. The question raised here is if criteria of the
second kind are necessary or even useful.

\section{Finding all invariants}

Consider a dimensional matrix with integer coefficients:
\[
\begin{array}{cc}
\begin{array}{c}
\phantom{D_{0}}\\
D_{1}\\
\cdots\\
D_{m}
\end{array} & \left\Vert \begin{array}{ccc}
q_{1} & \cdots & q_{n}\\
a_{11} & \cdots & a_{1n}\\
\cdots &  & \cdots\\
a_{m1} & \cdots & a_{mn}
\end{array}\right\Vert \end{array}.
\]
Recall \cite{dan} that finding an invariant product $q_{1}^{c_{1}}\ldots q_{n}^{c_{n}}$
of $q_{1},\dots,q_{n}$ is equivalent to finding a vector $\left(c_{1},\ldots,c_{n}\right)\in\mathbb{Z}^{n}$
such that
\begin{equation}
c_{1}\mathbf{v}_{1}+\ldots+c_{n}\mathbf{v_{n}}=\mathbf{0},\label{eq:inv}
\end{equation}
where $\mathbf{v}_{i}$ is the column vector $\left[a_{1i},\ldots,a_{mi}\right]^{\mathrm{T}}$
corresponding to $q_{i}$ and $\mathbf{0}$ is the column vector with
$n$ zeros. 

Let $\mathbf{D}$ be the matrix whose columns are the column vectors
$\mathbf{v}_{1},\ldots,\mathbf{v}_{n}$, let $r$ be the rank of $\mathbf{D}$,
and assume without loss of generality that the last $r$ columns of
$\mathbf{D}$ are independent as column vectors. Recall from linear
algebra that the solution space for (\ref{eq:inv}), $\mathbb{R}_{0}^{n}\!\left(\mathbf{D}\right)=\left\{ \left(c_{1},\ldots,c_{n}\right)\mid c_{i}\in\mathbb{R},\sum_{i}c_{i}\mathbf{v}_{i}=\mathbf{0}\right\} $,
has dimension $n-r$, and has a basis of $n-r$ vectors of the forms
\begin{equation}
\left(1,0,\ldots,0,b_{11},\ldots,b_{1r}\right),\dots,\left(0,\ldots,0,1,b_{\left(n-r\right)1},\ldots,b_{\left(n-r\right)r}\right).\label{eq:sol-inv}
\end{equation}
Thus, every solution of the equation system is a unique linear combination
of these $n$-tuples. As $\mathbf{D}$ is an integer matrix, all $b_{ij}$
are ratios of integers. Let $b_{i}>0$ be the lowest common denominator
of $b_{i1},\ldots,b_{ir}$. Then,
\begin{gather*}
I=\left\{ \left(b_{1},0,\ldots,0,c_{11},\ldots,c_{1r}\right),\dots,\left(0,\ldots,0,b_{n-r},c_{\left(n-r\right)1},\ldots,c_{\left(n-r\right)r}\right)\right\} ,
\end{gather*}
where $c_{ij}=b_{i}b_{ij}$, is a basis for $\mathbb{R}_{0}^{n}\!\left(\mathbf{D}\right)$
with only integer entries, and any integer solution of (\ref{eq:inv})
is a unique linear combination with integer coefficients of elements
of $I$. 

This means that every invariant corresponding to a solution of (\ref{eq:inv})
is a product of powers of $n-r$ invariants
\begin{gather*}
q_{1}^{b_{1}}q_{2}^{0}\ldots q_{n-r}^{0}q_{n-r+1}^{c_{11}}\ldots q_{n}^{c_{1r}},\dots,q_{1}^{0}\ldots q_{n-r-1}^{0}q_{n-r}^{b_{1}}q_{n-r+1}^{c_{\left(n-r\right)1}}\ldots q_{n}^{c_{\left(n-r\right)r}}.
\end{gather*}
 Since we can disregard factors of the form $q^{0}$ \cite{dan},
these invariants can be written in the form
\[
q_{i}^{b_{i}}\prod_{j=1}^{r}q_{n-r+j}^{c_{ij}}\qquad\left(i=1,\ldots,n-r\right),
\]
where $b_{i}>0$, and each invariant can be reduced further to an
invariant of the form
\begin{equation}
q_{i}^{b_{i}}\prod_{k=1}^{r_{i}}q_{n-r+j_{k}}^{c_{ij_{k}}}\qquad\left(i=1,\ldots,n-r,\quad0\leq r_{i}\leq r\right),\label{eq:pigroup}
\end{equation}
where $b_{i}>0,c_{ij_{k}}\neq0$ for all $i,j_{k}$, and $b_{i},c_{ij_{k}}$
are relatively prime.

For each such invariant, the column vectors corresponding to quantities
with non-zero exponents make up a set of non-independent column vectors.
Such a set is even a minimal set of non-independent column vectors;
the column vectors in any subset are independent because otherwise
all exponents $b_{i},c_{ij_{k}}$ would not be non-zero.

Thus, for (\emph{a}) every maximal (possibly empty) set of independent
quantities or column vectors, there is (\emph{b}) a (non-empty) set
of (non-empty) minimal sets of non-independent quantities or column
vectors, corresponding to invariants of the form (\ref{eq:pigroup}).
With terminology inspired by matroid theory, we can call a maximal
set of independent quantities or column vectors a \emph{basis set},
a minimal set of non-independent quantities or column vectors a \emph{circuit
set}, and we can rephrase the last sentence by saying that for any
basis set given by the dimensional matrix there are one or more corresponding
circuit sets. The set of circuit sets is in one-to-one correspondence
with (\emph{c}) a set of integer tuples $\left(c_{1},\ldots,c_{n}\right)$
called\emph{ circuit tuples} and (\emph{d}) a set of invariants $q_{1}^{c_{1}}\ldots q_{n}^{c_{n}}$
called \emph{circuit invariants}, where $c_{1},\ldots,c_{n}$ are
relatively prime, or, equivalently, (\emph{d'}) a set of corresponding
reduced circuit invariants of the form (\ref{eq:pigroup}). Note that
the sets of repeating quantities discussed in Section 1 are precisely
the basis sets.

Hence, if we are given a dimensional matrix $\mathbf{D}$, and we
form the union over all basis sets of the sets of circuit invariants
corresponding to sets of circuit sets associated with the current
basis set, we are sure to include all invariants used in dimensional
analysis based on $\mathbf{D}$. 

Consider, for example, the main example in Section 2. By inspection
of the dimensional matrix, we find ten basis sets, and there are two
circuit invariants of the form (\ref{eq:pigroup}) for each basis
set. The basis sets and corresponding circuit invariants are shown
below.

\begin{eqnarray*}
Basis\; set &  & Circuit\; invariants\\
\left\{ \nicefrac{\triangle P}{\ell},\rho,\mu\right\}  &  & d^{3}\left(\nicefrac{\triangle P}{\ell}\right)\rho\mu^{-2},u^{3}\left(\nicefrac{\triangle P}{\ell}\right)^{-1}\rho^{2}\mu^{-1}\\
\left\{ \nicefrac{\triangle P}{\ell},\rho,d\right\}  &  & \mu^{2}\left(\nicefrac{\triangle P}{\ell}\right)^{-1}\rho^{-1}d^{-3},u^{2}\left(\nicefrac{\triangle P}{\ell}\right)^{-1}\rho d^{-1}\\
\left\{ \nicefrac{\triangle P}{\ell},\mu,d\right\}  &  & \rho\left(\nicefrac{\triangle P}{\ell}\right)\mu^{-2}d^{3},u\left(\nicefrac{\triangle P}{\ell}\right)^{-1}\mu d^{-2}\\
\left\{ \nicefrac{\triangle P}{\ell},\rho,u\right\}  &  & \mu\left(\nicefrac{\triangle P}{\ell}\right)\rho^{-2}u^{-3},d\left(\nicefrac{\triangle P}{\ell}\right)\rho^{-1}u^{-2}\\
\left\{ \nicefrac{\triangle P}{\ell},\mu,u\right\}  &  & \rho^{2}\left(\nicefrac{\triangle P}{\ell}\right)^{-1}\mu^{-1}u^{3},d^{2}\left(\nicefrac{\triangle P}{\ell}\right)\mu^{-1}u^{-1}\\
\left\{ \nicefrac{\triangle P}{\ell},d,u\right\}  &  & \rho\left(\nicefrac{\triangle P}{\ell}\right)^{-1}d^{-1}u^{2},\mu\left(\nicefrac{\triangle P}{\ell}\right)^{-1}d^{-2}u\\
\left\{ \rho,\mu,d\right\}  &  & \left(\nicefrac{\triangle P}{\ell}\right)\rho\mu^{-2}d^{3},u\rho\mu^{-1}d\\
\left\{ \rho,\mu,u\right\}  &  & \left(\nicefrac{\triangle P}{\ell}\right)\rho^{-2}\mu u^{-3},d\rho\mu^{-1}u\\
\left\{ \rho,d,u\right\}  &  & \left(\nicefrac{\triangle P}{\ell}\right)\rho^{-1}du^{-2},\mu\rho^{-1}d^{-1}u^{-1}\\
\left\{ \mu,d,u\right\}  &  & \left(\nicefrac{\triangle P}{\ell}\right)\mu^{-1}d^{2}u^{-1},\rho\mu^{-1}du
\end{eqnarray*}
 The union of the ten sets of invariants is
\[
\left\{ \frac{\left(\nicefrac{\triangle P}{\ell}\right)\rho d^{3}}{\mu^{2}},\frac{\left(\nicefrac{\triangle P}{\ell}\right)\mu}{\rho^{2}u^{3}},\frac{\left(\nicefrac{\triangle P}{\ell}\right)d}{\rho u^{2}},\frac{\left(\nicefrac{\triangle P}{\ell}\right)d^{2}}{\mu u},\frac{\rho du}{\mu}\right\} ^{\pm1}.
\]
The set constructed in this way\emph{ }is a \emph{sufficient set of
minimal invariants} for $\mathbf{D}$, called the \emph{unified basis}
for $\mathbf{D}$, denoted $\mathcal{U\!}\left(\mathbf{D}\right)$.

Alternatively, one can start with the set of circuit sets for $\mathbf{D}$.
This set can be shown to be:
\[
\left\{ \left\{ \nicefrac{\triangle P}{\ell},\rho,\mu,d\right\} ,\left\{ \nicefrac{\triangle P}{\ell},\rho,\mu,u\right\} ,\left\{ \nicefrac{\triangle P}{\ell},\rho,d,u\right\} ,\left\{ \nicefrac{\triangle P}{\ell},\mu,d,u\right\} ,\left\{ \rho,\mu,d,u\right\} \right\} .
\]
Using $\mathbf{D}$, we find the set of pairs of circuit invariants
corresponding to this set of circuit sets:
\[
\left\{ \frac{\left(\nicefrac{\triangle P}{\ell}\right)\rho d^{3}}{\mu^{2}},\frac{\left(\nicefrac{\triangle P}{\ell}\right)\mu}{\rho^{2}u^{3}},\frac{\left(\nicefrac{\triangle P}{\ell}\right)d}{\rho u^{2}},\frac{\left(\nicefrac{\triangle P}{\ell}\right)d^{2}}{\mu u},\frac{\rho du}{\mu}\right\} ^{\pm1}.
\]
We call a set of invariants obtained in this way the \emph{circuit
basis} for $\mathbf{D}$, denoted $\mathcal{C\!\left(\mathbf{D}\right)}$.
The Cocoa script in the Appendix implements an algorithm that can
be used to calculate $\mathcal{C\!\left(\mathbf{D}\right)}$ (and
hence all circuit sets and circuit tuples as well).

The unified basis $\mathcal{U\!}\left(\mathbf{D}\right)$ is obviously
a subset of $\mathcal{C\!\left(\mathbf{D}\right)}$. Equality does
not always hold; for example, the dimensional matrix
\[
\begin{array}{cc}
\begin{array}{c}
\\
X
\end{array} & \left\Vert \begin{array}{cc}
q_{1} & q_{2}\\
1 & -1
\end{array}\right\Vert \end{array}
\]
has the unified basis $\left\{ q_{1}q_{2}\right\} $ but the circuit
basis $\left\{ q_{1}q_{2}\right\} ^{\pm1}$. In the case considered
here $\mathcal{U\!\left(\mathbf{D}\right)}=\mathcal{C\!\left(\mathbf{D}\right)}$,
however, and this is the typical situation.

How big can $\mathcal{C}\!\left(\mathbf{D}\right)$ be for a dimensional
matrix of rank $r$ with $n$ columns? It can be shown that this number
attains its maximum when the quantities in all sets of $r$ quantities
are independent. Then all circuit sets contain exactly $r+1$ quantities,
the number of circuit sets is $\binom{n}{r+1}$, and the number of
circuit invariants is $2\binom{n}{r+1}$. On the other hand, there
are $n-r$ active minimal invariants in any relation of the form (\ref{eq:Phi-rel}),
and all these invariants are circuit invariants from different pairs
of minimal invariants, so the number $\left\Vert \mathcal{C}\left(\mathbf{D}\right)\right\Vert $
of \emph{pairs} of circuit invariants satisfies the inequalities

\[
n-r\leq\left\Vert \mathcal{C}\!\left(\mathbf{D}\right)\right\Vert \leq\binom{n}{r+1}.
\]
For the four dimensional matrices in Sections \ref{sub:1-1}, 4 and
\ref{sec:Sets-of-representations}, we have

\[
\begin{array}{cccccc}
Section & n & r & n-r & \left\Vert \mathcal{C}\!\left(\mathbf{D}\right)\right\Vert  & \binom{n}{r+1}\\
---- & - & - & --- & --- & ---\\
2\:\left(1\right) & 5 & 3 & 2 & 5 & 5\\
2\:\left(2\right) & 4 & 3 & 1 & 1 & 1\\
4 & 5 & 2 & 3 & 8 & 10\\
5 & 5 & 3 & 2 & 3 & 5
\end{array}
\]

\subsection*{Remark 1}

Let $\sqsubseteq$ be a partial order on $\mathbb{Z}_{0}^{n}\!\left(\mathbf{D}\right)-\left\{ \mathbf{0}\right\} $,
the set of non-trivial integer solutions of (\ref{eq:inv}), such
that $\left(x_{1},\ldots,x_{n}\right)\sqsubseteq\left(y_{1},\ldots,y_{n}\right)$
if and only if $0\leq x_{i}\leq y_{i}$ or $0\geq x_{i}\geq y_{i}$
for $i=1,\ldots,n$. (Equivalently, $x_{i}y_{i}\geq0$ and $\left|x_{i}\right|\leq\left|y_{i}\right|$
for $i=1,\ldots,n$.) The Graver basis for $\mathbf{D}$, denoted
$\mathcal{G}\!\left(\mathbf{D}\right)$, is the set of all minimal
elements of $\mathbb{Z}_{0}^{n}\!\left(\mathbf{D}\right)-\left\{ \mathbf{0}\right\} $
under this partial order.

Corresponding to the circuit basis of invariants $\mathcal{C}\!\left(\mathbf{D}\right)$
as defined above, there is a \emph{circuit basis of tuples} $\mathcal{C}_{T}\!\left(\mathbf{D}\right)$
such that $\left(c_{1},\ldots,c_{n}\right)\in\mathcal{C}_{T}\!\left(\mathbf{D}\right)$
if and only if $q_{1}^{c_{1}}\ldots q_{n}^{c_{n}}\in\mathcal{C}\!\left(\mathbf{D}\right)$.
Every circuit tuple in $\mathcal{C}_{T}\!\left(\mathbf{D}\right)$
is a clearly a minimal element of $\mathbb{Z}_{0}^{n}\!\left(\mathbf{D}\right)-\left\{ \mathbf{0}\right\} $
under $\sqsubseteq$, so $\mathcal{C}_{T}\!\left(\mathbf{D}\right)\subset\mathcal{G}\!\left(\mathbf{D}\right)$.
Equality does not hold, however. For example, the matrix $\mathbf{D}=\left[\begin{array}{ccc}
1 & 2 & 1\end{array}\right]$ has the Graver basis 
\[
\left\{ \pm\left(2,-1,0\right),\pm\left(1,0,-1\right),\pm\left(0,1,-2\right),\pm\left(1,-1,1\right)\right\} ,
\]
but $\left(1,-1,1\right)$ and $\left(-1,1,-1\right)$ are not circuit
tuples.

In a recent e-print \cite{ather}, sets of invariants for dimensional
matrices corresponding to Graver bases for these matrices are presented.
It is noted that ``the Graver basis gives a full set of ... primitive
invariants'' for a dimensional matrix (p. 10). It has been shown
here, however, that it suffices to consider circuit bases.

\subsection*{Remark 2}

The terms 'unified basis' and 'circuit basis' are actually somewhat
misleading, since the invariants in such bases are not independent
in the usual sense. For example, $\rho^{2}\mu^{-1}u^{3}=\rho d^{-1}u^{2}\cdot\rho ud\mu^{-1}$,
and using such dependencies, the four representations of the relation
$\nicefrac{\triangle P}{\ell}=\Psi\!\left(\rho,\mu,d,u\right)$ can
be derived from each other by simple transformations. For example,
\begin{gather*}
\frac{\nicefrac{\triangle P}{\ell}}{\rho^{2}\mu^{-1}u^{3}}=\,\phi_{2}\!\left(\rho ud\mu^{-1}\right)\Longleftrightarrow\frac{\nicefrac{\triangle P}{\ell}}{\rho^{2}\mu^{-1}u^{3}}\rho ud\mu^{-1}=\,\phi_{2}\!\left(\rho ud\mu^{-1}\right)\rho ud\mu^{-1}\\
\Longleftrightarrow\frac{\nicefrac{\triangle P}{\ell}}{\rho d^{-1}u^{2}}=\,\phi_{3}\!\left(\rho ud\mu^{-1}\right).
\end{gather*}

This manifest equivalence of representations, which is not surprising
in view of the fact that they represent the same functional relation,
has been invoked in a sophisticated justification of the tradition
of limiting dimensional analysis to one relation and $n-r$ invariants.
The basic argument is that since all representations of the functional
relation are equivalent, it suffices to consider one of them, corresponding
to one, arbitrarily chosen, basis set (\cite{sonin}, p. 48). (A similar
argument can be found already in Buckingham's original article on
the $\Pi$ theorem \cite{bucking}, p. 362.) The examples in Section
2 should suffice to cast some doubt on this argument, however, and
it is further weakened by examples of dimensional analysis presented
in Sections 4 and 5.

In more abstract terms, we can regard the bijection $\left(c_{1},\ldots,c_{n}\right)\mapsto q_{1}^{c_{1}}\ldots q_{n}^{c_{n}}$
as an isomorphism between the circuit basis of tuples $\mathcal{C}_{T}\!\left(\mathbf{D}\right)$
and the circuit basis $\mathcal{C}\!\left(\mathbf{D}\right)$. Since
$\mathcal{C}_{T}\!\left(\mathbf{D}\right)$ is basically an $\left(n-r\right)$-dimensional
solution space (over $\mathbb{Z}$), $\mathcal{C}\!\left(\mathbf{D}\right)$
is an $\left(n-r\right)$-dimensional space as well with operations
defined by\linebreak{}
 $q_{1}^{c_{1}}\ldots q_{n}^{c_{n}}\cdot q_{1}^{d_{1}}\ldots q_{n}^{d_{n}}=q_{1}^{c_{1}+d_{1}}\ldots q_{n}^{c_{n}+d_{n}}$
and $\left(q_{1}^{c_{1}}\ldots q_{n}^{c_{n}}\right)^{a}=q_{1}^{ac_{1}}\ldots q_{n}^{ac_{n}}$,
and any choice of a set of repeating quantities is equivalent to a
choice of a basis for $\mathcal{C}\!\left(\mathbf{D}\right)$, where
the $n-r$ basis elements correspond to the $n-r$ invariants in a
representation of the given functional relation. While any two bases
are equivalent in the sense that the elements of one can be expressed
in terms of elements of the other, this does not mean that there cannot
be any benefit from considering more than one basis. Similarly, the
equivalence of representations in the sense indicated does not mean
that only one representation of the functional relation should be
considered.

\section{\label{sec2-}Using some of the invariants}

Consider quantities $q_{1},...,q_{n}$, a dimensional matrix $\mathbf{D}$
for these quantities, and a functional relation $q_{1}=\Psi\!\left(q_{2},...,q_{n}\right)$.
The invariants, obtained by dimensional analysis, which appear in
representations of this relation are all contained in $\mathcal{C}\!\left(\mathbf{D}\right)$,
but all invariants in $\mathcal{C}\!\left(\mathbf{D}\right)$ are
not used in the representations. This is basically because if a dependent
quantity has been designated, we should disregard basis sets which
include this quantity, because the column corresponding to the dependent
quantity is always linearly dependent on the columns corresponding
to the quantities in a basis set \cite{dan}. Hence, only invariants
associated with the non-disregarded basis sets will appear in the
representations of the functional relation.

Consider, for example, the ten basis sets and associated circuit invariants
in the preceding section. If $\nicefrac{\triangle P}{\ell}$ is designated
as the dependent quantity, meaning that we assume a functional relation
$\nicefrac{\triangle P}{\ell}=\Psi\!\left(\rho,\mu,d,u\right)$ to
hold, we should disregard the first six sets basis sets with associated
sets of invariants. For the reader's convenience, the four last basis
sets with associated sets of invariants are reproduced here:
\begin{eqnarray*}
Basis\; set &  & Circuit\; invariants\\
\left\{ \rho,\mu,d\right\}  &  & \left(\nicefrac{\triangle P}{\ell}\right)\rho\mu^{-2}d^{3},u\rho\mu^{-1}d\\
\left\{ \rho,\mu,u\right\}  &  & \left(\nicefrac{\triangle P}{\ell}\right)\rho^{-2}\mu u^{-3},d\rho\mu^{-1}u\\
\left\{ \rho,d,u\right\}  &  & \left(\nicefrac{\triangle P}{\ell}\right)\rho^{-1}du^{-2},\mu\rho^{-1}d^{-1}u^{-1}\\
\left\{ \mu,d,u\right\}  &  & \left(\nicefrac{\triangle P}{\ell}\right)\mu^{-1}d^{2}u^{-1},\rho\mu^{-1}du
\end{eqnarray*}

For each of these four basis sets there is a representation of $\nicefrac{\triangle P}{\ell}=\Psi\!\left(\rho,\mu,d,u\right)$
involving the associated invariants; for $\left\{ \rho,\mu,d\right\} $
we have $\nicefrac{\triangle P}{\ell}=\frac{\mu^{2}}{\rho d^{3}}\,\phi_{1}\!\left(\frac{\rho ud}{\mu}\right)$,
and so forth. These are the same formulas as given in Section 2, except
that we obtain $\nicefrac{\triangle P}{\ell}=\frac{\rho u^{2}}{d}\,\phi_{3}\!\left(\frac{\mu}{\rho ud}\right)$
instead of $\nicefrac{\triangle P}{\ell}=\frac{\rho u^{2}}{d}\,\phi_{3}\!\left(\frac{\rho ud}{\mu}\right)$
.

Thus, the basic reason why the circuit basis includes more invariants
than are needed in the representations of any single functional relation
is that a circuit basis can accommodate all possible functional relations
among the quantities considered, or equivalently, all choices of a
dependent quantity.

The following example, originally constructed by White and Lewalle,
is adapted from \cite{white}. The displacement $S\!\left(t\right)$
of a falling body as a function of elapsed time $t$ is given by the
differential equation $S''\!\left(t\right)=g$, where $g$ is the
(local) constant of gravity. This equation has the solution $S\!\left(t\right)=S_{0}+V_{0}t+\frac{1}{2}gt^{2},$
where $S_{0}=S\!\left(0\right)$ and $V_{0}=S'\!\left(0\right)$.
Thus, $S\!\left(t\right)=\Psi\!\left(S_{0},V_{0},g\right)\!\left(t\right)$.
It is instructive to perform a dimensional analysis as if $\Psi$
were an unknown function. The dimensional matrix is
\[
\begin{array}{cc}
\begin{array}{c}
\\
L\\
T
\end{array} & \left\Vert \begin{array}{ccccc}
S\!\left(t\right) & S_{0} & V_{0} & g & t\\
1 & 1 & 1 & 1 & 0\\
0 & 0 & -1 & -2 & 1
\end{array}\right\Vert .\end{array}
\]
This matrix has rank $2$, so $n-r=3$, but the circuit basis has
16 elements, namely $S\!\left(t\right)S_{0}^{-1}$, $S\!\left(t\right)V_{0}^{-2}g$,
$S\!\left(t\right)V_{0}^{-1}g^{-1}$, $S\!\left(t\right)g^{-1}t^{-2}$,
$S_{0}V_{0}^{-2}g$, $S_{0}V_{0}^{-1}t^{-1}$, $S_{0}g^{-1}t^{-2}$,
$V_{0}g^{-1}t^{-1}$ and their inverses, and there are $9$ basis
sets: $\left\{ S\!\left(t\right),V_{0}\right\} $, $\left\{ S\!\left(t\right),g\right\} $,
$\left\{ S\!\left(t\right),t,\right\} $, $\left\{ S_{0},V_{0}\right\} $,
$\left\{ S_{0},g\right\} $, $\left\{ S_{0},t,\right\} $, $\left\{ V_{0},g\right\} $,
$\left\{ V_{0},t\right\} $, $\left\{ g,t\right\} $. 

As usual, we disregard basis sets containing the dependent quantity
$S\!\left(t\right)$, and in this case we also disregard basis sets
containing $t$ because what we want to find out is how the functional
relation $t\mapsto S\!\left(t\right)$ depends on the parameters $S_{0}$,
$V_{0}$ and $g$. The remaining basis sets and invariants are the
following
\begin{eqnarray*}
Basis\; set &  & Circuit\; invariants\\
\left\{ S_{0},V_{0}\right\}  &  & S\!\left(t\right)S_{0}^{-1},gS_{0}V_{0}^{-2},tV_{0}S_{0}^{-1}\\
\left\{ S_{0},g\right\}  &  & S\!\left(t\right)S_{0}^{-1},V_{0}^{2}S_{0}^{-1}g^{-1},t^{2}S_{0}^{-1}g\\
\left\{ V_{0},g\right\}  &  & S\!\left(t\right)V_{0}^{-2}g,S_{0}V_{0}^{-2}g,tV_{0}^{-1}g
\end{eqnarray*}
 There are three corresponding representations of the functional relation
$S\!\left(t\right)=\Psi\!\left(S_{0},V_{0},g\right)\!\left(t\right)$:
\[
\begin{cases}
S\left(t\right)=S_{0}\,\Phi_{1}\!\left(\frac{g}{S_{0}^{-1}V_{0}^{2}},\frac{t}{S_{0}V_{0}^{-1}}\right) & \left(\mathfrak{e}\right)\\
S\left(t\right)=S_{0}\,\Phi_{2}\!\left(\frac{V_{0}}{\sqrt{S_{0}g}},\frac{t}{\sqrt{S_{0}g^{-1}}}\right) & \left(\mathfrak{f}\right)\\
S\left(t\right)=V_{0}^{2}g^{-1}\,\Phi_{3}\!\left(\frac{S_{0}}{V_{0}^{2}g^{-1}},\frac{t}{V_{0}g^{-1}}\right) & \left(\mathfrak{g}\right)
\end{cases}.
\]

It is pointed out in \cite{white} that while $\left(\mathfrak{e}\right)$,
$\left(\mathfrak{f}\right)$ and $\left(\mathfrak{g}\right)$ are
representations of the same functional relation, and thus can be said
to contain the same information, this information is presented in
different ways, making it possible to draw different kinds of conclusions
from the representations. For example, a plot of $S\!\left(t\right)/S_{0}$
as a function of $t/\left(S_{0}V_{0}^{-1}\right)$ for different values
of $g/\left(S_{0}^{-1}V_{0}^{2}\right)$ shows the effect of $g/\left(S_{0}^{-1}V_{0}^{2}\right)$
on the functional relation $t/\left(S_{0}V_{0}^{-1}\right)\mapsto S\!\left(t\right)/S_{0}$.
Hence, this plot shows the effect of $g$ on the functional relation
$t\mapsto S\!\left(S_{0},V_{0},g\right)\!\left(t\right)$ for constant
values of $S_{0}$ and $V_{0}$. Plots corresponding to $\left(\mathfrak{f}\right)$
and $\left(\mathfrak{g}\right)$ similarly show the effect of $V_{0}$
and $S_{0}$, respectively, on this functional relation.

\subsection*{Remark}

This example is one of many showing that one can get more informative
results from dimensional analysis by considering more than one functional
relation and more than $n-r$ invariants. Yet, this fact is not reflected
by the step-by-step algorithm for dimensional analysis presented in
\cite{white}. This algorithm only describes how to obtain one relation
involving $n-r$ invariants, so there is a tension between traditional
principles and advanced practice in this textbook's exposition of
dimensional analysis. It is argued here that this tension should be
eliminated by modifying the principles.

\section{\label{sec:Sets-of-representations}Sets of representations of functional
relations as equation systems}

We have seen examples of how the existence of more than $n-r$ minimal
invariants allows alternative representations of a relation of the
form $q=\Psi\!\left(q_{1},\ldots,q_{n}\right)$ to be derived by means
of dimensional analysis, but alternative representations are not only
alternatives. These alternative representations can be regarded as
an \emph{equation system}, from which we can obtain more information
about $\Psi$ than is available from the representations considered
separately. It may even be possible to determine $\Psi$ (up to a
multiplicative constant) from this equation system and available additional
information. This will be shown by means of an example, adapted from
\cite{bridgma}.

Let two bodies with mass $m_{1}$ and $m_{2}$ revolve around each
other in circular orbits under influence of their mutual gravitational
attraction. Let $d$ denote their distance and $t$ the time of revolution.
We want to derive a relation which shows how $t$ depends on relevant
parameters.

Preliminary considerations indicate that we should also include the
universal gravitational constant $G$ among the parameters, so that
we obtain the dimensional matrix
\[
\begin{array}{cc}
\begin{array}{c}
\\
L\\
T\\
M
\end{array} & \left\Vert \begin{array}{ccccc}
t & d & m_{1} & m_{2} & G\\
0 & 1 & 0 & 0 & 3\\
1 & 0 & 0 & 0 & -2\\
0 & 0 & 1 & 1 & -1
\end{array}\right\Vert \end{array}.
\]
The circuit basis is
\[
\left\{ \frac{m_{1}}{m_{2}},\:\frac{t^{2}}{d^{3}m_{1}^{-1}G^{-1}},\:\frac{t^{2}}{d^{3}m_{2}^{-1}G^{-1}}\right\} ^{\pm1},
\]
the basis sets corresponding to the functional relation $t=\Psi\!\left(d,m_{1},m_{2},G\right)$
are $\left\{ d,m_{1},G\right\} $ and $\left\{ d,m_{2},G\right\} $,
and the two sets of invariants associated with these two basis sets
are $\left\{ t^{2}m_{1}d^{-3}G,m_{2}m_{1}^{-1}\right\} $ and $\left\{ t^{2}m_{2}d^{-3}G,m_{1}m_{2}^{-1}\right\} $,
respectively, so the corresponding equation system is 
\[
\begin{cases}
t^{2}=\frac{d^{3}}{m_{1}G}\,\Phi_{1}\!\left(\frac{m_{2}}{m_{1}}\right) & \left(\mathfrak{h}\right)\\
t^{2}=\frac{d^{3}}{m_{2}G}\,\Phi_{2}\!\left(\frac{m_{1}}{m_{2}}\right) & \left(\mathfrak{h}'\right)
\end{cases}.
\]
This shows that $t^{2}$ is proportional to $d^{3}$ and inversely
proportional to $G$, so we have derived Kepler's third law in a special
case.

There is more information hidden in this equation system, however.
In view of the symmetry between the two revolving bodies we may assume
that $\Phi_{1}=\Phi_{2}=\Phi$. Multiplying $\left(\mathfrak{h}\right)$
and $\left(\mathfrak{h}'\right)$ by $\left(m_{1}G\right)/d^{3}$
and setting $x=m_{1}/m_{2}$, we obtain the functional equation
\[
\Phi\!\left(1/x\right)=x\,\Phi\!\left(x\right),
\]
 which has solutions of the form
\[
\Phi\!\left(x\right)=\frac{K}{1+x}.
\]
 Substituting this in $\left(\mathfrak{h}\right)$ or $\left(\mathfrak{h}'\right)$,
we get
\[
t^{2}=\frac{Kd^{3}}{G\left(m_{1}+m_{2}\right)}\qquad\mathrm{or}\qquad t=C\sqrt{\frac{d^{3}}{G\left(m_{1}+m_{2}\right)}}.
\]
There are more examples in \cite{dan} showing that dimensional analysis
can lead further than generally recognized if all relevant invariants
are taken into account, which strengthens the conclusion that in general
we should consider more than $n-r$ invariants and more than one relation
of the form (\ref{eq:Phi-rel}) in dimensional analysis.

\rule[0.5ex]{1\columnwidth}{1pt}

\lyxaddress{\textsc{\small{}Dan Jonsson, Department of Sociology and Work Science,
University of Gothenburg, SE 405 30 Gothenburg, Sweden.}}
\begin{lyxcode}
\newpage{}
\end{lyxcode}

\part*{Appendix}

Shown below is a simple CoCoA-5%
\footnote{\textbf{John Abbott, Anna Maria Bigatti, Giovanni Lagorio.} CoCoA-5:
a system for doing Computations in Commutative Algebra. Available
at http://cocoa.dima.unige.it.%
} script to help calculate the circuit basis. (Edit to specify the
dimensional matrix!\texttt{)}

\texttt{$\vphantom{x}$}~\\
\texttt{\textquotedbl{}\textbackslash{}n\textquotedbl{}; Use QQ; }~\\
\texttt{-{}- Edit list of quantities: }~\\
\texttt{vList := {[} \textquotedbl{}S(t)\textquotedbl{},\textquotedbl{}S0\textquotedbl{},\textquotedbl{}V0\textquotedbl{},\textquotedbl{}g\textquotedbl{},\textquotedbl{}t\textquotedbl{}
{]}; }~\\
\texttt{-{}- Edit list of columns in dimensional matrix: }~\\
\texttt{cList := {[} {[}1,0{]},{[}1,0{]},{[}1,-1{]},{[}1,-2{]},{[}0,1{]}
{]}; }~\\
\texttt{n := len( vList ); }~\\
\texttt{r := rank( matrix( cList ) ); }~\\
\texttt{For i := 1 To r+1 Do }~\\
\texttt{$\hphantom{x}\hphantom{x}$lsv := subsets( vList, i); }~\\
\texttt{$\hphantom{x\hphantom{x}}$lsm := subsets( cList, i ); }~\\
\texttt{$\hphantom{x\hphantom{x}}$m := len( lsm ); }~\\
\texttt{$\hphantom{x\hphantom{x}}$For j := 1 To m Do }~\\
\texttt{$\hphantom{x\hphantom{x\hphantom{x\hphantom{x}}}}$sm := lsm{[}
j {]}; }~\\
\texttt{$\hphantom{x\hphantom{x\hphantom{x\hphantom{x}}}}$sm := Transposed(
matrix( sm ) ); }~\\
\texttt{$\hphantom{x\hphantom{x\hphantom{x\hphantom{x}}}}$If rank(
sm ) = i-1 Then }~\\
\texttt{$\hphantom{x\hphantom{x\hphantom{x\hphantom{x}}}\hphantom{x\hphantom{x}}}$L
:= LinKerBasis( sm ); }~\\
\texttt{$\hphantom{x\hphantom{x\hphantom{x\hphantom{x\hphantom{x\hphantom{x}}}}}}$If
count( L{[} 1 {]} , 0 ) = 0 Then }~\\
\texttt{$\hphantom{x\hphantom{x\hphantom{x\hphantom{x\hphantom{x\hphantom{x}}}\hphantom{x\hphantom{x}}}}}$lut
:= {[} {]}; }~\\
\texttt{$\hphantom{x\hphantom{x\hphantom{x\hphantom{\hphantom{x\hphantom{x\hphantom{x}}}\hphantom{x\hphantom{x}}}}}}$For
k := 1 To len( lsv{[} j {]} ) Do }~\\
\texttt{$\hphantom{x\hphantom{x\hphantom{x\hphantom{x\hphantom{x\hphantom{x\hphantom{x\hphantom{x}}}\phantom{x}\hphantom{x}}}}}}$append(
ref lut, lsv{[} j,k {]} ); }~\\
\texttt{$\hphantom{x}\hphantom{x}\hphantom{x}\hphantom{x}\hphantom{x}\hphantom{x}\hphantom{x}\hphantom{x}\hphantom{x}\hphantom{x}$append(
ref lut, -L{[} 1,k {]} ); }~\\
\texttt{$\hphantom{x}\hphantom{x}\hphantom{x}\hphantom{x}\hphantom{x}\hphantom{x}\hphantom{x\hphantom{x}}$EndFor;
}~\\
\texttt{$\hphantom{x}\hphantom{x}\hphantom{x}\hphantom{x\hphantom{x}\hphantom{x}\hphantom{x}\hphantom{x}}$lut;
}~\\
\texttt{$\hphantom{x}\hphantom{x}\hphantom{x}\hphantom{x\hphantom{x}\hphantom{x}}$EndIf;
}~\\
\texttt{$\hphantom{x\hphantom{x}\hphantom{x}\hphantom{x}}$EndIf;
}~\\
\texttt{$\hphantom{x}\hphantom{x}$EndFor; }~\\
\texttt{EndFor;}~\\
$\vphantom{x}$\\
Output from this script (corresponding to $2\cdot8$ circuit invariants):\texttt{}~\\
\texttt{$\vphantom{x}$}\\
\texttt{{[}\textquotedbl{}S(t)\textquotedbl{}, -1, \textquotedbl{}S0\textquotedbl{},
1{]} }~\\
\texttt{{[}\textquotedbl{}S(t)\textquotedbl{}, 1, \textquotedbl{}V0\textquotedbl{},
-2, \textquotedbl{}g\textquotedbl{}, 1{]} }~\\
\texttt{{[}\textquotedbl{}S(t)\textquotedbl{}, -1, \textquotedbl{}V0\textquotedbl{},
1, \textquotedbl{}t\textquotedbl{}, 1{]} }~\\
\texttt{{[}\textquotedbl{}S(t)\textquotedbl{}, -1/2, \textquotedbl{}g\textquotedbl{},
1/2, \textquotedbl{}t\textquotedbl{}, 1{]} }~\\
\texttt{{[}\textquotedbl{}S0\textquotedbl{}, 1, \textquotedbl{}V0\textquotedbl{},
-2, \textquotedbl{}g\textquotedbl{}, 1{]} }~\\
\texttt{{[}\textquotedbl{}S0\textquotedbl{}, -1, \textquotedbl{}V0\textquotedbl{},
1, \textquotedbl{}t\textquotedbl{}, 1{]} }~\\
\texttt{{[}\textquotedbl{}S0\textquotedbl{}, -1/2, \textquotedbl{}g\textquotedbl{},
1/2, \textquotedbl{}t\textquotedbl{}, 1{]} }~\\
\texttt{{[}\textquotedbl{}V0\textquotedbl{}, -1, \textquotedbl{}g\textquotedbl{},
1, \textquotedbl{}t\textquotedbl{}, 1{]}}
\end{document}